\documentstyle{amsppt}

\magnification\magstep1
\NoRunningHeads
\hsize6.9truein
\pageheight{23 truecm}
\input btxmac.tex               
\bibliographystyle{hplain}
\loadbold
\input epsf
\input diagrams



\def\ideal#1.{I_{#1}}
\def\ring#1.{\Cal O_{#1}}
\def\proj#1.{\Bbb P(#1)}
\def\pr #1.{\Bbb P^{#1}}
\def\af #1.{\Bbb A^{#1}}
\def\Hz #1.{\Bbb F_{#1}}
\def\Hbz #1.{\overline{\Bbb F}_{#1}}
\def\fb#1.{\underset #1 \to \times}
\def\ten#1.{\underset #1 \to \otimes}
\def\res#1.{\underset {\ \ring #1.} \to \otimes}
\def\au#1.{\operatorname {Aut}\,(#1)}
\def\deg#1.{\operatorname {deg } (#1)}
\def\pic#1.{\operatorname {Pic}\,(#1)}
\def\pico#1.{\operatorname{Pic}^0(#1)}
\def\picg#1.{\operatorname {Pic}^G(#1)}
\def\ner#1.{NS (#1)}
\def\rdown#1.{\llcorner#1\lrcorner}
\def\rup#1.{\ulcorner#1\urcorner}
\def\cone#1.{\operatorname {NE}(#1)}
\def\ccone#1.{\overline{\operatorname {NE}}(#1)}
\def\coef#1.{\frac{(#1-1)}{#1}}
\def\vit#1.{D_{\langle #1 \rangle}}
\def\mm#1.{\overline {M}_{0,#1}}
\def\H1#1.{H^1(#1,{\ring #1.})}
\def\ac#1.{\overline {\Bbb F}_{#1}}
\def\mgn#1.#2.{\overline {M}_{#1,#2}}
\def\ilist#1.{{#1}_1,{#1}_2,\dots}

\def\adj#1.{\frac {#1-1}{#1}}
\def\spn#1.{\overline{#1}}
\def\pek#1.#2.{\Cal P^{#1}(#2)}
\def\plk#1.#2.{\Cal P^{\leq #1}(#2)}


\def\list#1.#2.{{#1}_1,{#1}_2,\dots,{#1}_{#2}}
\def\omitlist#1.#2.{{#1}_1,{#1}_2,\dots,\hat {{#1}_i}, \dots, {#1}_{#2}}
\def\omitlist0#1.#2.{{#1}_0,{#1}_1,\dots,\hat {{#1}_i}, \dots, {#1}_{#2}}
\def\loc#1.#2.{\Cal O_{#1,#2}}
\def\fderiv#1.#2.{\frac {\partial #1}{\partial #2}}
\def\map#1.#2.{#1 \longrightarrow #2}
\def\rmap#1.#2.{#1 \dasharrow #2}
\def\emb#1.#2.{#1 \hookrightarrow #2}
\def\non#1.#2.{\text {Spec }#1[\epsilon]/(\epsilon)^{#2}}
\def\Hi#1.#2.{\text {Hilb}^{#1}(#2)}
\def\sym#1.#2.{\operatorname {Sym}^{#1}(#2)}
\def\Hb#1.#2.{\text {Hilb}_{#1}(#2)}
\def\Hm#1.#2.{\Hom_{#1}(#2)}
\def\prd#1.#2.{{#1}_1\cdot {#1}_2\cdots {#1}_{#2}}
\def\Bl #1.#2.{\operatorname {Bl}_{#1}#2}
\def\pl #1.#2.{#1^{\otimes #2}}


\def\alist#1.#2.#3.{#1_1 #2 #1_2 #2\dots #2 #1_{#3}}
\def\zlist#1.#2.#3.{#1_0 #2 #1_1 #2\dots #2 #1_{#3}}
\def\lmap#1.#2.#3.{#1 \overset #2\to \longrightarrow #3}
\def\ses#1.#2.#3.{0\longrightarrow #1 \longrightarrow #2 \longrightarrow #3 
\longrightarrow 0}
\def\les#1.#2.#3.{0\longrightarrow #1 \longrightarrow #2 \longrightarrow #3}
\def\es#1.#2.#3.{#1 \longrightarrow #2 \longrightarrow #3}
\def\Hi#1.#2.#3.{\text {Hilb}^{#1}_{#2}(#3)}


\def\Hom{\operatorname{Hom}}

\def\deg{\operatorname{deg}}


\def\e{\Cal E}

\def\e1{E_1}
\def\e2{E_2}


\def\mapdown#1{\big\downarrow\rlap{$\vcenter
{\hbox{$\scriptstyle#1$}}$}}

\def\mapse#1{
{\vcenter{\hbox{$\mathop{\smash{\raise1pt\hbox{$\diagdown$}\!\lower7pt
\hbox{$\searrow$}}\vphantom{p}}\limits_{#1}\vphantom{\mapdown{}}$}}}}


\def\VR#1.{height#1pt&\omit&&\omit&&\omit&&\omit&&\omit&\cr}

\def\VRT#1.{height#1pt&\omit&&\omit&\cr}

\NoBlackBoxes
\topmatter
\title
Threefold thresholds
\endtitle
  
\author
 James M\raise 2.2pt\hbox{\text{\tenrm c}}Kernan and Yuri Prokhorov
\endauthor  
\address
James M\raise 2.2pt \hbox{\text {\sevenrm c}}Kernan: Department of Mathematics,
University of California at Santa Barbara,
Santa Barbara, CA 93106
\endaddress
\email
mckernan\@math.ucsb.edu
\endemail

\address
Yuri Prokhorov: Department of Algebra, Faculty of Mathematics,
Moscow State Lomonosov University, Leninskie Gory, GSP-2, Moscow,
119 992, Russia
\endaddress

\email
prokhoro\@mech.math.msu.su
\endemail
\abstract 
We prove that the set of accumulation points of thresholds in dimension three
is equal to the set of thresholds in dimension two, excluding one. 
\endabstract
\endtopmatter
\document

\head \S 1 Introduction\endhead

 We prove  

\proclaim{1.1 Theorem} The set of accumulation points of the log canonical 
threshold in dimension three is equal to the set of log canonical thresholds in dimension
two (excluding the number one).  
\endproclaim

 In fact we prove much more, we prove a similar Theorem in all dimensions, assuming the
MMP and a conjecture of Alexeev-Borisov, see \S 3.

  It is natural to attach numerical invariants to a singular variety $D$ which measure the
complexity of the singularities.  For example we can associate the multiplicity to a
hypersurface singularity.  Unfortunately this is rather a coarse invariant; for example
consider the infinite family $y^2+x^b=0$, $b\geq 2$, of plane curves of multiplicity two.
Clearly one needs to take account of higher order blow ups to distinguish these
singularities.

 One such invariant is the log canonical threshold.  Suppose that $D$ is embedded as a
hypersurface in the smooth variety $X$.  The log canonical threshold $c$ is then the
largest value of $c$ such that $K_X+cD$ is log canonical.  A similar definition pertains
when $X$ is log canonical and $D$ is an integral $\Bbb Q$-Cartier divisor.  If one picks a
log resolution $\map Y.X.$ and one writes
$$
K_Y+E=\pi^*K_X+\sum a_iE_i \qquad\qquad \tilde D+\sum b_iE_i=\pi^*D,
$$
for appropriate rational numbers $\list a.k.$ and $\list b.k.$ where $E$ is the
sum over all exceptional divisors $\list E.k.$ and $\tilde D$ is the strict
transform of $D$, then
$$
c=\min_i \frac {a_i}{b_i}. 
$$

 If $D$ is given as $y^2+x^b=0$, then the log canonical threshold $c=1/2+1/b$.  More
generally, if $D$ is a hypersurface of dimension $n$ and multiplicity $a\geq n$, then the
log canonical threshold $1/a\leq c\leq (n+1)/a$, so that the log canonical threshold is a
subtle refinement of the multiplicity. 

 Shokurov conjectured that the set of log canonical thresholds in dimension $n$ satisfies
ACC (the ascending chain condition).  For example, if $D$ is the plane curve $y^a+x^b=0$
then the log canonical threshold is $1/a+1/b$ and it is proved in \cite{Igusa77} (see also
\cite{Kuwata99a}) that the log canonical threshold of any irreducible
plane curve has this form.

 Now given any set that satisfies ACC, it is natural to ask for the set of accumulation
points.  Koll\'ar conjectured that the set of accumulation points in dimension $n$ is equal
to the set of accumulation points in dimension $n-1$.  We prove this conjecture in
dimension three.

 In fact the log canonical threshold is a very natural invariant; when $D$ is a
hypersurface, then the log canonical threshold is both the complex singular index of $D$
(a number that measures the asymptotic behaviour of the integrals of vanishing cycles) and
(minus) the largest root of the (reduced) Bernstein-Sato polynomial of $D$, see
\cite{Kollar95} for more details.   

 Besides being a very natural invariant of study when $D$ is a hypersurface in $X$, the
log canonical threshold plays an important role in inductive approaches to higher
dimensional geometry.  For example, one of the most important remaining issues in Mori's
program is to establish the existence and termination of flips.  After Shokurov's
groundbreaking work on the problem of existence of four-fold flips, some of the focus of
interest has turned to the problem of termination of flips; the log canonical threshold
and other closely related invariants should play a crucial role in this problem.

 Another reason for being interested in the log canonical threshold is its connection with
questions of birational geometry and rigidity.  Along these lines, \cite{Cheltsov01} and
\cite{CP02} prove some results and state some conjectures concerning the log canonical
thresholds of hypersurfaces in $\pr n.$.  In fact \cite{Pukhlikov02} claims that every
smooth hypersurface of degree $n\geq 4$ in $\pr n.$ is birationally superrigid and part of
the proof of this fact relies on the log canonical threshold.  Finally \cite{EM01},
\cite{EFM02a} and \cite{EFM02b} prove some of the conjectures in \cite{Cheltsov01} and
\cite{CP02}, and thereby establish some new cases of birational superrigidity. 

   In yet another direction, the log canonical threshold should be closely related to
recent work of Hassett, \cite{Hassett02}.  In this paper Hassett constructs a moduli space
of pairs, consisting of a curve together with a boundary, where the coefficients of the
boundary are rational.  A very natural question is to ask what happens as one varies the
coefficients.  As the coefficients are increased, some pairs transition from stable to
unstable, and it is necessary to blow up the corresponding loci.  In fact a pair becomes
unstables as it becomes not log canonical and the critical values, at which a transition
occurs, are therefore log canonical thresholds.

 Furthermore the set of all log canonical thresholds seems to have a rich structure of its
own.  The log canonical threshold has received some attention for these reasons.  Several
papers establish special cases of (1.1); \cite{Kollar92b} proves that the
largest log canonical threefold threshold is $41/42$; \cite{Kuwata99} identifies the set
of thresholds $c\geq 5/6$ for a surface $D$ in $\Bbb C^3$; similarly \cite{Prokhorov01a}
establishes that $5/6$ is the largest accumulation point for all pairs and
\cite{Prokhorov02} identifies all the accumulation points greater than $1/2$.

\cite{Kollar92b} proceeds by direct computation and Mori's explicit classification 
of terminal $3$-fold singularities.  \cite{Kuwata99a} proceeds by analysing all possible
exceptional divisors on a smooth surface that have log discrepancy zero with respect
to some divisor $K_X+\Delta$.  \cite{Kuwata99} proceeds by an explicit analysis of a
defining equation for $D$.  \cite{EM01}, \cite{EFM02a} and \cite{EFM02b} prove their
results by using jet schemes, a method introduced by
Mircea~Musta\c{t}\v{a}. \cite{Prokhorov01a} and \cite{Prokhorov02} use a method first
introduced by Alexeev \cite{Alexeev94}, who proved that the set of all log canonical
threefold thresholds satisfies ACC.  This paper builds on the ideas established in these
two papers, so that many of the ideas from this paper can be traced back to the work of
Alexeev.

Before we give precise statements of our results, we need to give some definitions.  As
these definitions are rather involved, it will probably help the reader to give some of
the ideas of the proof of (1.1).  As observed above, many of the ideas for the
argument below is due to Alexeev.  However what did come as a surprise to both authors is
that his argument also lends itself so well to identifying the set of accumulation points
of the log canonical threshold, a fact which does not seem to have been exploited before.

{\bf Thanks}: Much of the work for this paper was completed whilst the first author was
attending the programme in Higher Dimensional Geometry, at the Isaac Newton Institute.
Furthermore collaboration on this paper began as a result of the second author's talk at
this programme.  Both authors would like to thank the organisers of the programme for
creating such a stimulating and hospitable atmosphere to work in. The second author was
partially supported by the grant INTAS-OPEN-2000-269.  Part of this work was also
completed whilst the second author was visiting RIMS; the second author would like to
thank RIMS for their generous hospitality and support during his stay there.

\head \S 2 Sketch of Proof of (1.1)\endhead

The definition of log canonical states that the log discrepancy of various exceptional
divisors are all nonnegative.  Thus in the definition of the log canonical threshold, at
least one exceptional divisor has log discrepancy exactly equal to zero.  A natural way to
proceed then, would be to find a birational morphism $\pi\:\map Y.X.$ that extracts an
exceptional divisor $E$ of log discrepancy zero.  For example, if $D$ is the curve
$y^a+x^b=0$ then $\map Y.X.$ is precisely the weighted blow up of $X=\Bbb C^2$, with
weights $(a,b)$, in the given coordinates $x$, $y$.  In general, the existence of $\pi$
is guaranteed by the MMP (minimal model program), see (6.1).  

The next step is to restrict to the exceptional divisor $E$ and apply adjunction (here we
use the fact that $E$ has log discrepancy zero).  In this case the singular points of 
$Y$ along $E$ make a contribution.  For example, when $D=y^a+x^b=0$, $E$ is a copy of $\pr 1.$
and 
$$
(K_Y+E)|_E=K_E+\left (\adj a.\right ) p+\left (\adj b.\right ) q,
$$
where $p$ corresponds to the point of index $a$ and $q$ to the point of index $b$.
Similar formulae pertain in the general case; see \S 4 for more details.  If
$\tilde D$ denotes the strict transform of $D$, then $\tilde D\cdot E=1$.  Moreover by
definition $(K_Y+c\tilde D+E)\cdot E=0$ so that
$$
-2+\adj a.+\adj b.+c=0,
$$
and so
$$
c=\frac 1a+\frac 1b, 
$$
as claimed.  

 In general, then, when we restrict to $E$ we naturally get a log canonical pair
$(E,\Delta)$ where $K_E+\Delta$ is numerically trivial and the number $c$ appears in one
of the coefficients of $\Delta$.  Moreover, using the MMP, see \S 5, we can actually
reverse this process.  More often than not, it suffices to take the cone over the pair
$(E,\Delta)$.

 In this way, the MMP in dimension $n$ reduces the problem of computing the log canonical
threshold to the problem of working with pairs $(X,\Delta)$ of dimension $n-1$, such that
$K_X+\Delta$ is log canonical and numerically trivial.

 So consider the problem of trying to identify the set of accumulation points for the
coefficients of such pairs $(X,\Delta)$.  One reason that the dimension one case is easy
is that there is only one possibility for $X$, $X$ must be isomorphic to $\pr 1.$.  In
higher dimensions, it is not hard, running the MMP again, to reduce to the case where $X$
has Picard number one, so that at least $X$ is a Fano variety.  

  Now suppose that we had a sequence $(X_i,\Delta_i)$ of such pairs.  If $X_i$ is constant
then it is easy to see that the set of coefficients of $\Delta$ have to decrease and it is
easy to identify the limit of $c_i$.  A similar argument works, if the set $\{X_i\}$ forms a
bounded family.

 Thus we may assume that the set $\{X_i\}$ forms an unbounded family.  At this point we
invoke a Theorem of Alexeev, which is conjectured to hold in all dimensions
(3.7).  By his Theorem, it follows that the log discrepancy of $X_i$ must
approach zero.

 To see how to use this fact, suppose that in fact the log discrepancy were zero and not
just approaching zero.  That is suppose that we have a log canonical pair $(X,\Delta)$
which is numerically trivial but not kawamata log terminal.  In this case if we extracted
a divisor of log discrepancy zero and applied adjunction then we would have decreased the
dimension by one as before.

 Thus our aim is now clear.  Generate a component of $\Delta$ with coefficient one, so
that we can restrict to this component and apply adjunction.

 In the more general case, we extract a divisor of smallest log discrepancy, so that there
is a component of $\Delta$ whose coefficient is approaching one.  If we repeat this
process twice, then we get two divisors whose coefficient is approaching zero.  By taking
an appropriate combination of these two divisors, we generate our component of coefficient
one and we can then apply induction.  In practice there are some technical issues that
have to be dealt with to make this argument work.

\head \S 3 Precise definitions and Statement of results\endhead

  We work over an algebraically closed field of characteristic zero.   

Hopefully it will have become clearer the level of generality that we need to work at to
make our argument go through.  

 First, we are interested in two sets of numbers, those for thresholds and those when
$K_X+\Delta$ is numerically trivial.

 Second we have to allow $\Delta$ to have unusual coefficients, since even if we start
with one integral component of coefficient $c$, when we restrict to a component of
coefficient one and apply adjunction, we obtain more than one component, with various
coefficients.

 Third, we have to allow for $c$ to be a sequence of coefficients, rather than just one
coefficient.  Indeed we need two divisors whose coefficient is approaching one.

 Fourth, after running the MMP, we might lose information on the sequence of coefficients,
as some components might well be contracted.  In particular sometimes we can only
determine the set of coefficients and not the full sequence.  

 Fifth, there are some subtleties to do with kawamata log terminal versus log canonical.
On the one hand kawamata log terminal singularities are much better behaved in general.
Thus we can expect sharper results in this case.  On the other hand, in terms of induction
we want to allow log canonical singularities, as we want to restrict to a component of
coefficient one.

\definition{3.1 Definition}  Let $I$ be any subset of $(0,1]$.  $I_+$ denotes the set 
of all finite sums of elements of $I$ less than one. 

 $\Cal C$ denotes the set of finite sequences of elements of $(0,1]$ whilst $\Cal C^+$
denotes the subset of increasing finite sequences of elements.  Note the obvious inclusion
$(0,1]\subset \Cal C^+$.  Define a partial order on $\Cal C$ by the rule:\newline 
$\list x.m.<\list y.n.$ iff the two sequences are not equal, $m\leq n$ and $x_i\leq y_i$, for all $i\leq m$.
In addition, in the case of $\Cal C^+$, we require strict inequality for at least one $i$.  
\enddefinition

 We record the following easy result.  

\proclaim{3.2 Lemma} Let $c^{(i)}$ be a sequence in $\Cal C^+$.  

 Then $c^{(i)}$ is a strictly increasing sequence iff there is a strictly increasing sequence 
of numbers $x_i$, such that for every $i$, there is an index $j$, with $x_i=c^{(i)}_j$ 
\endproclaim
\demo{Proof} Clear.  \qed\enddemo

 Let $c^{(i)}$ be the first $i$ terms of the sequence of reciprocals of the natural
numbers.  Then $c^{(i)}$ is a strictly increasing sequence in $\Cal C$.

  Recall that a log canonical centre is the image of an exceptional divisor of log
discrepancy at most zero.   

\definition{3.3 Definition} We say that the sequence $\list c.m.\in \Cal C$ is 
{\bf associated} to the boundary $\Delta$, with coefficients $\list d.n.$, if for every
$i$ we may write
$$
d_i=\adj m_i. +\frac {f_i}{m_i}+ \sum_j \frac {k_{ij}c_j}{m_i},
$$
for some $f_i\in I_+$ and natural numbers $\list m.n.$ and $k_{ij}$, where for all $j$,
there is an $i$ such that $k_{ij}\neq 0$.  We say this sequence is of {\bf adjunction
type} if in addition $m_i\neq 1$ for some $i$; otherwise we say $c$ is {\bf ordinary}.

 Consider the following four conditions that can be imposed on a log canonical pair $(X,\Delta)$.
\roster 
\item"(T)" $X$ is $\Bbb Q$-factorial and every component of $\Delta$ contains a log 
canonical centre $Z$ of $(X,\Delta)$.   
\item"(N)" $X$ is complete and $K_X+\Delta$ is numerically trivial.
\item"(K)" condition N holds, and $K_X+\Delta$ is kawamata log terminal. 
\item"(B)" condition N holds, and for every $j$ there is a big  $\Bbb Q$-Cartier divisor whose support is 
contained in those components $\Delta_i$ of $\Delta$ such that $k_{ij}\neq 0$.
\endroster 
\enddefinition

  We can use (3.3) to build an assortment of subsets of $\Cal C$.  

\definition{3.4 Definition} To every log canonical pair, $(X,\Delta)$, 
where $X$ is a variety of dimension $n$, we take any sequence $c$ associated to $\Delta$.
We use the four letters T, K, N, B corresponding to the four conditions, with subscript
$n$.  For condition T we require that $c$ is ordinary and for condition $B$ we take $X$ of
dimension at most $n$.  Calligraphic letters, $\Cal T$, $\Cal N$, $\Cal K$, $\Cal B$
indicate that we take arbitrary sequences; the letters $T$, $N$, $K$, $B$,
that $c$ is a single number.
\enddefinition

 Thus $T_3(1)$ denotes the set of threefold thresholds and, for example, $\Cal K_5(I)$
denotes the set of all sequences associated to kawamata log terminal pairs $(X,\Delta)$,
where $X$ is a complete variety of dimension five and $K_X+\Delta$ is numerically trivial.

 Note that (3.4) involves some choices.  In particular we only allow ordinary sets
for the threshold; in fact it turns out, see (6.3), that even if we allow sets
of adjunction type, we get the same set of coefficients.  Some authorities replace
condition B by the condition that $X$ has Picard number one, see for example
\cite{Kollaretal}.  In fact the corresponding sets are conjecturally equal, see
(5.1) and (5.2).  We choose condition B as it behaves better in terms of
inductive arguments.  We allow $X$ to be of any dimension less than $n$, for condition B,
so that the inclusion $\Cal B_{n-1}(I)\subset \Cal B_n(I)$ is automatic.  Similar
inclusions for N and K are clear.

 Shokurov made the remarkable conjectures that these sets satisfy ACC.  Koll\'ar speculated
that even more is true, that the set of accumulation points in dimension $n$ should be
closely related to the corresponding sets in dimension $n-1$.

\proclaim{3.5 Conjecture} Let $I$ be any subset of $(0,1]$, which contains $1$.   
\roster 
\item $B_n(I)=K_n(I)=N_n(I)$.  
\item $\Cal T_{n+1}(I)=\Cal B_n(I)$.  
\item $I$ satisfies DCC iff $\Cal N^+_n(I)$ satisfies ACC iff $\Cal B_n(I)$ satisfies ACC.  
\item Suppose that the only accumulation point of $I$ is one.  Then the set of accumulation points 
of $N_n(I)$ is equal to $N_{n-1}(I)-\{1\}$.
\item $I$ satisfies DCC iff $\Cal K_n(I)$ satisfies ACC.  
\endroster 
\endproclaim

\example{3.6 Example} Let $X=\pr 2.$ and $\Delta=L+c_1L_1+c_2L_2+\dots c_mL_m$, 
where $\sum c_i=2$ and the lines $L,\list L.m.$ are in general position.  If we blow up
any of the points of intersection of $L$ and $L_i$ the exceptional divisor has coefficient
$c_i$ in the log pullback $\Gamma$ of $\Delta$.  Repeating this process arbitrarily often
we may find a surface $S$ and a divisor $\Gamma$ where the set of coefficients of $\Gamma$
is equal to $\{\list c.m.\}$ but where each coefficient appears arbitrarily many times.

 In particular $\Cal N_2(1)$ does not satisfy ACC.  
\endexample

 Recall the following very interesting conjecture, which is a weakening of a conjecture
due independently to Alexeev and Borisov (the standard version of this conjecture places
no restriction on the Picard number):

\proclaim{3.7 Conjecture} Fix an integer $n$ and a positive real 
number $\epsilon$.  

 Then the family of all varieties $X$ such that 
\roster 
\item the dimension of $X$ is $n$,
\item the log discrepancy of $X$ is at least $\epsilon$,
\item $-K_X$ is ample, and
\item the Picard number of $X$ is one,
\endroster 
is bounded. 
\endproclaim

Alexeev \cite{Alexeev94} proved (3.7) in dimension two and A. Borisov and L. Borisov
\cite{BB92} proved (3.7) for toric varieties.  Kawamata \cite{Kawamata92a}
proved (3.7) for Fano threefolds with Picard number one and terminal
singularities and Koll\'ar, Miyaoka, Mori and Takagi \cite{KMMT00} proved that all Fano
threefolds with canonical singularities are bounded.  Just recently Alexeev and Brion
\cite{AB03} have produced a preprint in which they claim to prove (3.7) for 
spherical varieties.

 We will also need the following conjecture, which is a modification of a conjecture made
in \cite{Alexeev94} by Alexeev for surfaces,

\proclaim{3.8 Conjecture} Fix an integer $n$ and a positive real 
number $\epsilon$.  

 Then the family of all pairs $(X,\rup \Delta.)$ such that 
\roster 
\item the dimension of $X$ is $n$, 
\item the log discrepancy of the pair $(X,\Delta)$ is at least $\epsilon$, and the coefficients of 
$\Delta$ are at least $\epsilon$ and at most $1-\epsilon$, 
\item $-(K_X+\Delta)$ is nef, 
\item $X$ is rationally connected, 
\endroster 
is bounded. 
\endproclaim

 Note that (3.8) is a much stronger version of (3.7), as any
kawamata log terminal Fano variety is automatically rationally connected.  Despite the
fact that (3.8) seems quite a bold conjecture, we give two pieces of evidence
for this conjecture; we prove that (3.8) holds if either one assumes the MMP and
(3.7) and we replace (3.8.4) by the condition that $-K_X$ is big, see
(10.3), or we replace (3.8.4) by the condition that $X$ is toric,
see (11.2).

\proclaim{3.9 Theorem} Assume the MMP in dimension $n$.  Then (3.5.1) 
and (3.5.2)$_{n-1}$ hold.  Further assume (3.7) holds in dimension
$n$.  Then (3.5.3-4) hold.  Finally assume that (3.8) holds in
dimension $n$.  Then (3.5.5) holds.
\endproclaim

\proclaim{3.10 Corollary} (3.5) holds for $n\leq 2$.  In particular 
(1.1) holds.  Further the largest accumulation point of $T_3$ is $5/6$.
\endproclaim

 In fact the set of accumulation points in dimension three is completely determined, see
(4.6).  Note that both the MMP and (3.7) are known to hold for toric
varieties.  We are therefore able to prove (3.5) in the toric case.  First some 
convenient notation.  

 \proclaim{3.11 Definition} We say that the pair $(X,\Delta)$ has a 
{\bf toroidal resolution} if there is a log resolution $\pi:\map Y.X.$ of the pair $(X,\Delta)$ 
which is toroidal with respect to some pair $(X,D)$.  
\endproclaim

 Now we introduce some more notation which modify the definitions of condition $N$ and
$T$.  We say that a pair $(X,\Delta)$ satisfies condition ${}^tT$ if it satisfies
condition $T$ and in addition the pair $(X,\Delta)$ has a toroidal resolution.  We say
that a pair $(X,\Delta)$ satisfies condition ${}^tN$ if it satisfies condition $N$ and in
addition $X$ is toric and every exceptional valuation of log discrepancy zero is a toric
valuation.  With these subsets, we obtain a modified version of (3.5), which we
denote by ${}^t$(3.5).

\proclaim{3.12 Theorem} ${}^t$(3.5) holds. 
\endproclaim

 A very special case of (3.12) is the result that thresholds for
non-degenerate hypersurfaces satisfy ACC (a non-degenerate hypersurface automatically
admits a toroidal resolution).  In fact one inspiration for stating (3.12) is
the preprint of \cite{Ishii00}, where it is proved that the log canonical threshold of
non-degenerate hypersurfaces satisfies ACC.

 It is relatively straightforward to extend the proofs of most of the results given in this
paper to characteristic $p$, although some proofs become a little more delicate, without
their exposition adding anymore insight.  In particular (6.3) does not hold in
characteristic $p$.  For this reason we only state and prove our results in characteristic
zero.

\head \S 4 Arithmetic of Adjunction\endhead

 We put together some of the arithmetic underlying adjunction and running the MMP in this 
section.   We start with a Proposition that we prove after a couple of Lemmas.  

\proclaim{4.1 Proposition} Suppose that we have a log canonical pair $(X,\Delta)$.
Let $Y$ be a component of $\Delta$ of coefficient one and define $\Gamma$ on $Y$ by
adjunction,
$$
K_Y+\Gamma=(K_X+\Delta)|_Y.
$$

 If the sequence $\list c.k.\in \Cal C$ is associated to the pair $(X,\Delta)$ and every
component of $\Delta$ intersects $Y$ then it is also associated to the pair $(Y,\Gamma)$.
\endproclaim

 We start with a useful definition and a well-known Lemma.

\definition{4.2 Definition} Let $I\subset (0,1]$ and set  
$$
D(I)=\{\, a\leq 1\,|\, a=\adj m.+\frac fm,m\in \Bbb N,f\in I_+\,\}.
$$
\enddefinition

\proclaim{4.3 Lemma} Suppose that we have a log canonical pair $(X,\Delta)$.
Let $Y$ be a component of $\Delta$ of coefficient one and define $\Gamma$ on $Y$ by
adjunction,
$$
K_Y+\Gamma=(K_X+\Delta)|_Y.
$$
 
 If the coefficients of $\Delta$ belong to $I$ then the coefficients of $\Gamma$ belong to
$D(I)$.
\endproclaim
\demo{Proof} Immediate from (16.6) of \cite{Kollaretal}.  \qed\enddemo

\proclaim{4.4 Lemma} Let $I$ be any subset of $(0,1]$, which contains $1$. 
Then 
\roster 
\item $D(I)_+=D(I)$.  
\item $D(D(I))=D(I)$.  
\item $I$ satisfies DCC iff $D(I)$ satisfies DCC.  
\endroster 
\endproclaim
\demo{Proof} Suppose that $g\in D(I)_+$.  It follows that there are real numbers 
$f_i\in I_+$ such that
$$
g=\sum \adj n_i.+\frac {f_i}{n_i}.
$$
 Now all but one $n_i=1$ else $g\geq 1$ and there is nothing to prove.  Thus $g$ has the form
$$
g=\adj n.+\frac fn
$$
for some $f\in I_+$.  Hence (1).  

 Suppose that $h\in D(D(I))$.  Then there is a real number $g\in D(I)_+$ and an integer
$m$ such that 
$$
h=\adj m.+\frac gm.
$$

 As $g\in D(I)_+$, $g\in D(I)$ by (1).  Thus $g$ has the form
$$
g=\adj n.+\frac fn
$$
for some $f\in I_+$.  Hence
$$
\align
h=\adj m.+\frac gm         &=\adj m.+\frac {(\adj n.+\frac fn)}m \\ 
                           &=\frac {mn-n+n-1}{mn}+\frac f{mn}\\
                           &=\adj mn.+\frac f{mn}=\adj r.+\frac fr.
\endalign
$$
where $r=mn$.  But $f\in I_+$, so that $h\in D(I)$.  Hence (2).  

 As $I\subset D(I)$, it is clear that if $D(I)$ satisfies DCC then so does $I$.  Now
suppose that $I$ satisfies DCC.  Let $\ilist b.$ be a non-increasing sequence in $D(I)$.
By assumption there are integers $m_{ij}$ and elements $a_i\in I_+$ such that
$$
b_i=\sum_j \adj m_{ij}.+\frac {a_i}{m_{ij}}. 
$$
If $b_i=1$ then there is nothing to prove.  Thus we may assume that at most one $m_{ij}\neq 1$ so that
$$
b=\sum \adj m_i.+\frac {a_i}{m_i}. 
$$
Thus we may assume that $m_i$ is constant.  In this case the result is clear. \qed\enddemo

\demo{Proof of (4.1)} This follows easily from (4.3) and 
(4.4).  \qed\enddemo

 We end this section with a couple of useful general results.  

\proclaim{4.5 Lemma} Suppose that $I\subset \Bbb R$ satisfies DCC and $K\subset \Bbb R$ satisfies 
ACC.  

Then 
$$
\multline X=\{\, \list c.p.\in \Cal C\,|\, \exists\, m_i,\, k_{ij}\in \Bbb N,\, a\in I, \, f_i\in I^+ \\
\ a+\sum _i \left (\adj m_i.+\frac {f_i}{m_i}+\sum _j \frac {k_{ij}c_j}{m_i}\right ) =k\in K, \\
\text {where each term $\adj m_i.+\frac {f_i}{m_i}+\sum _j \frac {k_{ij}c_j}{m_i}$ of the sum is at most one}\,\}
\endmultline
$$
satisfies ACC. 
\endproclaim
\demo{Proof} Suppose we have a sequence of elements $x^{(p)}$ of $X$.  If any terms 
$$
\adj m_i^{(p)}.+\frac {f_i^{(p)}}{m_i^{(p)}}+\frac {k_{ij}^{(p)}c_j^{(p)}}{m_i^{(p)}},
$$
appearing in the sum on the left are strictly increasing then we may absorb them into
$a^{(p)}$.  In this case $a^{(p)}$ is strictly increasing and so at least one term of
$c_i^{(p)}$ must appear elsewhere in the sum.  Possibly passing to a subsequence, it
follows that we may assume that the $m_i^{(p)}$ are constant.  Replacing each element $k$
of $K$ by $m_ik$, we may assume that $m_i=1$.  Replacing $I$ by $I_+\cup \{0\}$ and $a$ by
$a+\sum f_i$, we may assume that each $f_i=0$.  Replacing $K$ by all numbers of the form
$k-a$, where $k\in K$ and $a\in I$ we may assume that $I=\{0\}$ so that $a=0$.  In this
case the result is clear. \qed\enddemo
 
 In general it seem hard to identify the sets defined in (3.4). 
The one exception to this is $N_1$.

\proclaim{4.6 Lemma} 
$$
\multline N_1=N_1(1)=\{\, a\in (0,1] \,|\, \exists m_i,k_i\in \Bbb N,\ \sum \adj m_i.+\frac {k_ia}{m_i}=2 \\
\text {where each term $\adj m_i.+\frac {k_ia}{m_i}$ of the sum is at most one}\,\}
\endmultline
$$

 In particular the largest element of $N_1$, other than one, is $5/6$.  
\endproclaim
\demo{Proof} Suppose $a\in N_1$.  Then there is a curve $C$, and a
divisor $\Delta$ such that $K_C+\Delta$ is numerically trivial, where the coefficients of $\Delta$ have
the form
$$
\adj m_i.+\frac {k_ic}{m_i}.
$$
Assuming that $D\neq 0$ it follows that $C\simeq \pr 1.$ and so the sum of the
coefficients is two. \qed\enddemo

\head \S 5 Reduction to Picard number one\endhead

 We commence the proof of (3.9) in this section.  The proof of (3.9) consists in
showing that (3.5.1-5) holds, under suitable circumstances.  We adopt the convention that 
(3.9.i) refers to proving that (3.5.i) holds, under the hypotheses stated in 
(3.9). 

 The first thing to prove is (3.9.1).  

\proclaim{5.1 Lemma} Assume that the MMP holds in dimension $n$. 

 Let $(X,\Delta)$ be a log canonical pair satisfying condition N, where $X$ is of dimension $n$. 
Suppose that we are given big $\Bbb Q$-Cartier divisors $\list D.k.$ supported on $\Delta$.  

 Then we may run a MMP $\rmap X.Y.$, where $Y$ is projective with $\Bbb Q$-factorial log
terminal singularities, and there is an extremal contraction $\pi\:\map Y.Z.$, which is
not birational.  Further if $\Gamma$ is the log pullback of $\Delta$ then the pair
$(Y,\Gamma)$ satisfies condition N.  Moreover the strict transforms $\list G.k.$ of the
divisors $\list D.k.$ are big $\Bbb Q$-Cartier divisors.

 In particular, in the definition of $\Cal B_n(I)$, we may assume 
that $X$ is projective, log terminal and $\Bbb Q$-factorial of Picard number one.  
\endproclaim
\demo{Proof} Set $\Theta=\Delta-\epsilon(\sum D_i)$, for $\epsilon$ sufficiently small.  
Let $f\:\map W.X.$ be a log terminal model of the pair $(X,\Theta)$ and let $\Gamma$ be
the log pullback of $\Delta$.  Let $E$ be an exceptional divisor of $\map W.X.$.  As
$(X,\Delta)$ is log canonical, the centre of $E$ is not contained in any of the divisors
$D_j$.  Thus $f^*D_j$ is equal to the strict transform $G_j$ of $D_j$.  As $D_j$ is big,
then so is $G_j$.

 Replacing $(X,\Delta)$ by $(W,\Gamma)$ we may therefore assume that $X$ is projective and
$\Bbb Q$-factorial.  We run the $K_X$-MMP.  It is clear that the divisors $\list D.k.$
remain big.  Hence the first statement. 

 Restricting to the general fibre of $\pi$, and applying (4.1) and
induction on the dimension, gives the second statement.  \qed\enddemo

\proclaim{5.2 Lemma} Assume the MMP in dimension $n$.  

 Let $(X,\Delta)$ be a log canonical pair satisfying condition N, where $X$ is of
dimension $n$.  Let $A$ be a component of $\Delta$.  Then we may run a MMP $\rmap X.Y.$,
where $Y$ is projective with $\Bbb Q$-factorial log terminal singularities, and there is
an extremal contraction $\pi\:\map Y.Z.$, which is not birational.  Further if $\Gamma$ is
the log pullback of $\Delta$ then the pair $(Y,\Gamma)$ satisfies condition N and the
strict transform of $A$ dominates $Z$.

 Moreover in the definition of $N_n(I)$ we may assume that $X$ is projective,
$\Bbb Q$-factorial of Picard number one and that $(X,\Delta)$ is kawamata log terminal.
\endproclaim 
\demo{Proof} Passing to a log terminal model we may assume that $X$ is projective
$\Bbb Q$-factorial.  We run the $(K_X+\Delta-\epsilon A)$-MMP.  With this choice of MMP, $A$ is
automatically positive on every ray that we contract.  In particular we never contract
$A$.  We end up with a Mori fibre space $\pi\:\map X.Z.$.  $A$ dominates $Z$, as $A$ is
positive on the corresponding ray.

 Suppose we have a log pair $(X,\Delta)$, as in the definition of $N_n(I)$.  Let $A$ be a
component of $\Delta$ whose coefficient is associated to $c$.  Suppose that we contract a
component $B$ of coefficient one, or that a component of coefficient one does not dominate
$Z$.  Then at some stage of the MMP, $A$ must intersect $B$.  Restricting to $B$ and
applying (4.1) we are done by induction.

 Otherwise restrict to the general fibre of $\pi$, and apply (4.1) and
induction on the dimension. \qed\enddemo

\demo{Proof of (3.9.1)}  Note that if $X$ is $\Bbb Q$-factorial, of Picard number one, 
then any effective $\Bbb Q$-divisor is big and $\Bbb Q$-Cartier.  Thus the result follows easily
from (5.1) and (5.2).  \qed\enddemo

\head \S 6 Some cyclic covers\endhead

The next thing to prove is (3.9.2), assuming the existence of the MMP in
dimension $n$.  (3.9.2) involves proving that two sets are equal, which can be
broken down into proving that each set is a subset of the other.  One inclusion is easy
and well known.  We begin with a proof of this case, as it sets the stage for the reverse
inclusion.

\proclaim{6.1 Lemma} $\Cal T_n(I)\subset \Cal B_{n-1}(I)$.  
\endproclaim
\demo{Proof} Suppose that the pair $(X,\Delta)$ satisfies condition $T$.  Let $\pi\:\map W.X.$ 
extract a divisor $E$ of log discrepancy zero and let $Y$ be the normalisation of a
general fibre of $\pi|_E\:\map E.\pi(E).$.  Let $\Theta$ be the log pullback of $\Delta$.
We may assume that $W$ is $\Bbb Q$-factorial and that $\pi$ has relative Picard number
one.  In particular the restriction of every component of $\Theta$ to $Y$ is ample, whence
big and $\Bbb Q$-Cartier.  Define $\Gamma$ on $Y$ by adjunction.  This gives us a log pair
$(Y,\Gamma)$ which satisfies condition B.  By (4.1) the coefficients of
$\Gamma$ are of the appropriate form.  \qed\enddemo

\proclaim{6.2 Lemma} Suppose that the pair $(X,\Delta)$ satisfies condition 
N, where $X$ is $\Bbb Q$-factorial of Picard number one.  

 Then there is a pair $(Y,\Gamma)$ satisfying condition T, where $Y$ is of dimension one
more than $X$ and the coefficients of $\Gamma$ are equal to the coefficients of $\Delta$. 

 In particular if $c\in \Cal C$ is associated to $\Delta$ then $c$ is associated to $\Gamma$. 
\endproclaim
\demo{Proof} Pick an embedding of $X$ into projective space such that $X$ is projectively normal
and let $(Y,\Gamma)$ be the cone over $(X,\Delta)$.  As $X$ is projectively normal, $Y$ is
normal.  Suppose that $p$ is the vertex of the cone.  Note that every component of
$\Gamma$ contains $p$ and that $Y$ is $\Bbb Q$-factorial as $X$ has Picard number one.  

 Let $\map W.Y.$ blow up $p$, so that $W$ is a $\pr 1.$-bundle over $X$ and the
exceptional divisor $E$ is a copy of $X$.  Let $\Theta$ be the strict transform of
$\Gamma$.  We may write
$$
K_W+\Theta+E=\pi^*(K_Y+\Gamma)+aE
$$
where $a$ is by definition the log discrepancy of $E$.  Now 
$$
(K_W+\Theta+E)|_E=K_X+\Delta,
$$
under the natural identification, so that $a=0$.  It follows that $E$ is a log
canonical place and $p$ is a log canonical centre.  On the other hand $K_W+\Theta+E$ is
log canonical by inversion of adjunction and so the pair $(Y,\Gamma)$ is log canonical.
Thus the pair $(Y,\Gamma)$ satisfies condition T.  \qed\enddemo

 Note that we cannot apply (6.2) directly, as $c$ might not be ordinary.  In
characteristic zero, we can get around this by taking a cyclic cover.

\proclaim{6.3 Lemma} Suppose that the pair $(X,\Delta)$ satisfies condition T
and that the sequence $c$ is associated to $\Delta$.

 Then we may find a pair $(Y,\Gamma)$ satisfying condition T, where $Y$ has dimension at
most the dimension of $X$, $c$ is associated to $\Gamma$ and $c$ is ordinary.  Further we
may assume that every component of $\Gamma$ is Cartier.
\endproclaim
\demo{Proof} Note that condition $T$ is local about the generic point of $Z$, the log 
canonical centre.  In particular we may assume that $X$ is affine.  Cutting by hyperplanes
in $X$ we may assume that $Z$ is a point.

 Suppose the coefficient of one component $A$ of $\Delta$ is 
$$
\adj m.+\frac fm.
$$
Let $\pi:\map Y.X.$ be a cyclic cover of order $m$ which is totally ramified over the
generic point of $A$, but which is otherwise \'etale in codimension one.  Let $\Gamma$ be
the log pullback of $\Delta$.  If $B$ is the inverse image of $A$, then $\pi^*A=mB$ so
that the coefficient of $B$ is equal to $f$.  Moreover $K_Y+\Gamma$ is log canonical and
$q$ the inverse image of $p$ is a log canonical centre of $K_Y+\Gamma$.  Taking another
cyclic cover of $Y$, which is \'etale in codimension one, we may assume that $B$ is
Cartier.

 We are therefore done by induction on the number of components of $\Gamma$. \qed\enddemo

\demo{Proof of (3.9.2)} By (6.1) it suffices to prove 
$\Cal B_{n-1}(I)\subset\Cal T_n(I)$.  Suppose we are given a pair $(X,\Delta)$ that
satisfies condition B, where $X$ has dimension $n-1$.  By (5.1) we may assume
that $X$ is projective, $\Bbb Q$-factorial, of Picard number one.  By (6.2) we
may find a pair $(Y,\Gamma)$ that satisfies condition T, where $Y$ has dimension $n$.  By
(6.3) we may assume that if $c$ is associated to the pair $(X,\Delta)$ then $c$
is associated to $(Y,\Gamma)$ and moreover $c$ is ordinary.  Thus $c\in\Cal
T_n(I)$. \qed\enddemo

(6.3) has another useful application.  

\proclaim{6.4 Lemma} Every element of $T_{n-1}(I)-\{1\}$ is an accumulation 
point of $T_n(I)$.  
\endproclaim
\demo{Proof} Suppose we have a pair $(X,\Delta)$ satisfying condition T.  By (6.3) 
we may assume that every component of $\Delta$ is Cartier.  The result now follows easily
from (8.21.2) of \cite{Kollar95}.  \qed\enddemo

\head \S 7 The exchange trick\endhead 

 The following simple trick will prove quite useful.  

\proclaim{7.1 Lemma} Let $(X_i,\Delta_i,A_i+B_i)$ be a sequence of triples, 
where $(X_i,\Delta_i)$ is a log pair satisfying condition N, $A_i$ and $B_i$ are Weil
divisors on $X_i$ and where each $X_i$ is $\Bbb Q$-factorial, of Picard number one.
Suppose that the coefficients of $A_i$ and $B_i$ converge.

 Then, possibly passing to a subsequence, there is a divisor $\Gamma_i$, with the following 
properties
\roster 
\item The pair $(X_i,\Gamma_i)$ also satisfies condition N. 
\item $\Delta_i-\Gamma_i$ is supported on $A_i+B_i$.  
\item the limits of the coefficients of $A_i$ and $B_i$ also converge to the same respective limits. 
\item Possibly switching $A_i$ and $B_i$, we may assume that either the coefficient of $A_i$ is constant,
or there is a log canonical centre contained in $A_i$.   
\item If the coefficient of $B_i$ in $\Delta_i$ is increasing (respectively decreasing) then the 
coefficient of $B_i$ in $\Gamma_i$ is increasing (respectively decreasing).
\item If the coefficients of both $A_i$ and $B_i$ are either strictly increasing or strictly 
decreasing then the coefficient of $B_i$ in $\Gamma_i$ is not constant.
\endroster 
\endproclaim
\demo{Proof} Since $X_i$ is $\Bbb Q$-factorial of Picard number one, $A_i=\lambda_iB_i$,
for some $\lambda_i$.  Possibly passing to a subsequence and switching $A_i$ and $B_i$ we
may assume that $\lambda_i|a-a_i|\leq |b-b_i|$, where $a_i$ (respectively $b_i$) is the
coefficient of $A_i$ (respectively $B_i$) and $a$ (respectively $b$) is the limit of $a_i$
(respectively $b_i$).  Set $\Theta_t=\Delta_i+t(a-a_i)A_i-t\lambda_i(a-a_i)B_i$.  When
$t=0$, we get $\Delta_i$.  If $K_X+\Theta_1$ is log canonical for $t=1$, then set
$\Gamma_i=\Theta_1$.  Otherwise pick the largest value of $t$ such that $K+\Theta_t$ is
log canonical and let $\Gamma_i=\Theta_t$.  The result is clear. \qed\enddemo

\proclaim{7.2 Lemma} Let $(X_i,\Delta_i)$ be a sequence of pairs, 
where $(X_i,\Delta_i)$ is a log pair satisfying condition N, and each $X_i$ is $\Bbb
Q$-factorial, of Picard number one.  Suppose that every coefficient of $\Delta_i$ is 
monotonic.  

Then, possibly passing to a tail of the sequence, there are divisors $\Gamma_i$ and $H_i$, with the
following properties
\roster 
\item The pair $(X_i,\Gamma_i)$ also satisfies condition N. 
\item $H_i$ belongs to a base point free linear system.    
\item The support of $\Gamma_i$ is contained in the support of $\Delta_i+H_i$.  
\item The coefficients of $\Gamma_i$ and $\Delta_i$ have the same limit.  
\item The only decreasing coefficient of $\Gamma_i$ is the coefficient of $H_i$.  
\endroster 
\endproclaim
\demo{Proof} As every coefficient of $\Delta_i$ has a limit, then, possibly passing to a 
tail of the sequence, we may form a divisor $\Theta_i$ on $X_i$, with support equal to the
support of $\Delta_i$, whose coefficients are equal to the limiting coefficients.  We may
write
$$
\Theta_i-\Delta_i=P_i-N_i,
$$
where $P_i$ and $N_i$ have disjoint support and both are effective.  Clearly the limits of
the coefficients of $P_i$ and $N_i$ are both zero, $\Theta_i-P_i$ is effective and the
components of $N_i$ correspond to those divisors whose coefficients are decreasing.  As
$X_i$ is $\Bbb Q$-factorial of Picard number one, $N_i$ is ample and so we may find
$m_i>0$, increasing, such that $m_iN_i$ is base point free.  Pick $H_i$ a general element
of the corresponding linear system.  Set
$$
\Gamma_i=(\Theta_i-P_i)+\frac 1{m_i}H_i.
$$
The result is then clear.  \qed\enddemo

\proclaim{7.3 Lemma} Suppose that $\Cal T^+_n(I)$ satisfies ACC.  

 Let $(X_i,\Delta_i)$ be a sequence of log pairs satisfying condition N, where each $X_i$
is $\Bbb Q$-factorial, of Picard number one and dimension $n$.  Suppose that the
coefficients of $\Delta_i$ form an increasing sequence of $\Cal C^+$.

 Then we may find a divisor $\Gamma_i$, such that $(X_i,\Gamma_i)$ also satisfies
condition N and exactly one coefficient of $\Gamma_i$ is increasing.  If further
$\Cal B^+_{n-1}(I)$ satisfies ACC then at most one coefficient of $\Delta_i$ is
approaching one.
\endproclaim
\demo{Proof} Suppose that there is more than one component of $\Delta_i$ whose coefficient 
is increasing.  We apply the exchange trick to two such components.  As the coefficients
of $\Gamma_i$ are increasing and we are assuming that $\Cal T^+_n(I)$ satisfies ACC, the
second possibility for (7.1.3) is ruled out.  Hence the first statement follows
from repeated application of (7.1) and (3.2).

 If there is more than one component whose coefficient is approaching one, then we may
assume that there is a component of coefficient one.  Restricting to this component and
applying (4.1) gives a contradiction. \qed\enddemo

\proclaim{7.4 Lemma} Assume the MMP in dimension $n$.  Suppose that 
both $\Cal B_n(I)$ and $\Cal T^+_n(I)$ satisfy ACC.

 Let $(X_i,\Delta_i)$ be a sequence of log pairs satisfying condition N, where each $X_i$
is $\Bbb Q$-factorial, of Picard number one and dimension $n$.  Suppose that the
coefficient $a_i$ of one component $A_i$ of $\Delta_i$ is approaching one.  Suppose that
$c_i$ approaches a limit $c$, where $c_i$ is a sequence associated to $\Delta_i$.

 Then $c$ belongs to $\Cal B_{n-1}(I)$.  
\endproclaim
\demo{Proof} As we are assuming that $\Cal B_n(I)$ satisfies ACC there is at 
least one component $B_i$ of $\Delta_i$ whose coefficient $b_i$ is decreasing.

 We first apply (7.3).  Thus we may assume that $\Delta_i$ has exactly one 
component $B_i$ whose coefficient is decreasing, that this coefficient is approaching zero, 
and that $B_i$ is chosen general in a base point free linear system.  

We apply (7.1) and replace $\Delta_i$ by the sequence of divisors $\Gamma_i$ so
produced.  Either the coefficient of $A_i$ is one or the coefficient of $B_i$ is zero.
The coefficient of $A_i$ must be one otherwise there is no component whose coefficient is
decreasing.  Hence we may assume that the coefficient of $A_i$ is one.  In this case we
may restrict to $A_i$ and apply (4.1). \qed\enddemo

\head \S 8 From Boundedness to ACC\endhead

\proclaim{8.1 Lemma}  Suppose we have a sequence $(X_i,\Delta_i)$ of 
log pairs and let $c_i$ be the sequence associated to $\Delta_i$.   

 If $\{X_i\}$ belongs to a bounded family, then $c_i$ is not an increasing sequence
of $\Cal C$. 

 If further each $X_i$ is $\Bbb Q$-factorial of dimension $n$ and Picard number one, $c_i$
is a number (rather than a sequence of numbers) and the only accumulation point of $I$ is
one, then $\lim c_i\in B_{n-1}(I)$.
\endproclaim
\demo{Proof} By assumption we may find a curve $C_i$ belonging to the smooth locus of $X_i$ 
such that $-K_{X_i}\cdot C_i\leq M$ independently of $i$.  As $C_i$ lies in the smooth locus of
$X_i$ the intersection number of $C_i$ with every component of $\Delta_i$ is a positive integer.
The possible set of values of $-K_{X_i}\cdot C_i$ is finite and any finite set satisfies ACC.
The first statement follows by (4.5). 

 Now suppose that the only accumulation point of $I$ is one and that $X_i$ has Picard
number one.  It follows that the coefficient of at least one component of $\Delta_i$ is
approaching one and the second statement follows by (7.4).  \qed\enddemo

\demo{Proof of (3.9.3)} Suppose not.  Then there is an increasing sequence $c_i$, either 
in $\Cal N^+_n(I)$ or $\Cal B_n(I)$.  Let $(X_i,\Delta_i)$ be the corresponding sequence
of log pairs.  We first reduce to the case of Picard number one.  There are two cases.  If
$c_i$ is an increasing sequence in $\Cal B_n(I)$, then by (5.1) we may assume
that $X_i$ has Picard number one.  Otherwise $c_i$ belongs to $\Cal N^+_n(I)$.  For each
$i$, by (3.2), we may find a component $A_i$ of $\Delta_i$ whose
coefficient is increasing.  By (5.1) we may assume that $X_i$ is $\Bbb
Q$-factorial, of Picard number one.

 As we are reduced to the case of Picard number one, then in fact we may assume that
$c_i\in\Cal B_n(I)$.  By (8.1) we may assume that $\{X_i\}$ does not form a
bounded set.  As we are assuming (3.7) it follows that the log discrepancy of
$X_i$ is approaching zero.  In particular we may assume that the log discrepancy is less
than one.  By (5.2) we reduce to the case where $\Delta_i$ contains a component
$A_i$ whose coefficient is approaching one.  

 Applying (8.1) a second time, we may extract $\pi_i\:\map Y_i.X_i.$ another
divisor $B_i$ whose log discrepancy is approaching zero.  By induction we may assume that
$\Cal T_n(I)$ satisfies ACC so that we may assume that
$K_{Y_i}+\Gamma_i+(1-a_i)A_i+(1-b_i)B_i$ is log canonical, where $\Gamma_i$ is the log
pullback of $\Delta$ and $a_i$ (respectively $b_i$) is the coefficient of $A_i$
(respectively $B_i$).

 We apply (5.2) to find a sequence of flips $\rmap Y_i.Y_i'.$ followed by an
extremal contraction $\pi'\:\map Y_i'.X_i'.$ which is not small, on which $B_i'$ is
positive.  We let $\Gamma_i'$ be the strict transform of $\Gamma_i$ and $\Delta_i'$ the
pushforward of $\Delta_i$.  

 If $\pi_i'$ is not birational then we restrict to the general fibre and apply
(4.1) and induction on the dimension.  If $\pi'$ is birational and $A_i$ is
not contracted then we are reduced to the case where we have two components whose
coefficient is approaching one and we apply (7.3) to derive a contradiction.

 Otherwise we apply (9.2).  Possibly switching the roles of $A_i$ and $B_i$ we may
assume that $K_{Y_i}+\Theta_i$ is log canonical and $\pi$-trivial, where $A_i$ has
coefficient one in $\Theta_i$ and the coefficient of $B_i$ is still increasing.
Restricting to $A_i$ and applying (4.1) and induction on the dimension we
are done.  \qed\enddemo

\demo{Proof of (3.9.4)} Suppose that we have a sequence $c_i\in N_n(I)$.  
Let $(X_i,\Delta_i)$ be the corresponding sequence of log pairs.  By (5.2) we may
assume that $X_i$ is $\Bbb Q$-factorial of Picard number one and that $(X_i,\Delta_i)$ is
kawamata log terminal.

By (8.1) we may assume that $\{X_i\}$ does not form a bounded set.  As we are
assuming (3.7) it follows that the log discrepancy of $X_i$ is approaching
zero.  In particular we may assume that the log discrepancy is less than one.  Let
$\pi_i:\map Y_i.X_i.$ extract a divisor $A_i$ of minimal log discrepancy.  

 By (5.2) we reduce to the case where $\Delta_i$ contains a component whose
coefficient is approaching one.  The only thing to check is that in this process we have
not contracted every component associated to $c_i$.  But (3.9.3) implies that
there must be a component whose coefficient is decreasing, and so there must be a
component whose coefficient involves $c_i$.  Now apply (7.4).\qed\enddemo

\head \S 9 Two ray game\endhead

We employ the following simple variants of the two ray game.  We introduce some convenient
notation.  Let $Y$ be a normal projective variety with $\Bbb Q$-factorial singularities
and Picard number two.  Then the cone of curves of $Y$ is spanned by two rays $R$ and $S$.
Now suppose that we are given a sequence of flips $g\:\rmap Y.Y'.$.  Then the cone of
curves of $Y'$ is also spanned by two rays $R'$ and $S'$.  Possibly switching the roles of
$R$ and $S$, we may assume that $R$ is spanned by the first curve that is flipped in $Y$
and that $S'$ is spanned by the last flipping curve in $Y'$.

\proclaim{9.1 Lemma} Let $D$ be any divisor in $Y$, which is pseudo-effective.
Suppose that $D$ is negative on $S$.  Then $D$ is positive on $R$, and $D'$ is positive on
$R'$ and negative on $S'$, where $D'$ denotes the strict transform of $D$.
\endproclaim
\demo{Proof} As $D$ is negative on $S$, it is certainly not numerically trivial over $Y$.  
As $D$ is pseudo-effective, but not numerically trivial, it is positive on a covering
family of curves.  As the cone of curves is spanned by $R$ and $S$, $D$ must be positive
on $R$.  By induction on the number of flips, it suffices to consider the case when $g$ is
a flip.  As $D$ is positive on $R$ and $g$ is the flip of $R$, $D'$ is negative on $S'$.
By what we have already proved it follows that $D'$ is positive on $R'$. \qed\enddemo

 We introduce some more notation.  Suppose in addition that we have a boundary $\Delta$ so
that we have a log canonical pair $(Y,\Delta)$, where $K_Y+\Delta$ is numerically trivial.
Assume that both rays $S$ and $R'$ are contractible.  Let $f\:\map Y.X.$ be the
contraction of $S$ and $f'\:\map {Y'}.{X'}.$ the contraction of $R'$.  Assume that $f$ and
$f'$ are both birational and not small.  Let $E$ and $E'$ be the corresponding exceptional
divisors.  We may write $K_Y+\Delta=K_Y+\Gamma+aE+a'E'$, where $a$ is the coefficient of
$E$ and $a'$ is the coefficient of $E'$.

\proclaim{9.2 Lemma} Suppose that $(Y,\Delta)$ is kawamata log terminal. 

 Then either $K_Y+\Gamma+E+b'E'$ is trivial on $S$, for some $a'<b'<1$, or
$K_{Y'}+\Gamma+bE+E'$ is trivial on $R'$, for some $a<b<1$.
\endproclaim
\demo{Proof} Consider the divisor $D=K_Y+\Gamma+E+E'$.  Suppose that $D$ is $S$-positive.  
As $E$ is negative on $S$, it follows that $E'$ is $S$-positive.  Thus $K_Y+\Gamma+E+b'E'$
is trivial on $S$ for some $b'<1$.  On the other hand
$K_Y+\Gamma+E+a'E'=K_Y+\Delta+(1-a)E$ is negative on $S$.  Thus $a'<b'$.

 Thus we may assume that $D$ is $S$-negative.  On the other hand as $D$ is numerically
equivalent to an effective linear combination of $E$ and $E'$, $D$ is certainly
pseudo-effective.  Thus by (9.1) $D'$ is positive on $R'$.  Applying the
same argument as above, it follows that $K_{Y'}+\Gamma+bE+E'$ is trivial on $R'$, for some
$a<b<1$. \qed\enddemo

\head \S 10 Proof of (3.9.5)\endhead

 To prove (3.9.5) we need the following result.  

\proclaim{10.1 Theorem} Let $X$ be a projective variety.   

 There is a rational map $\pi\:\rmap X.B.$ with the following properties.  
Let $U$ be the open subset of $X$ where $\pi$ is defined.  Then $\pi|_U$ is proper and 
the fibres of $\pi|_U$ are rationally chain connected.  Moreover $B$ is not uniruled.  
\endproclaim
\demo{Proof} Let $\pi$ be the maximally rationally connected chain fibration, whose existence 
is established in \cite{KMM92}, see also \cite{Kollar96}.  By the main Theorem of \cite{GHS03},  
$B$ is not uniruled. \qed\enddemo

(10.1) has the following very useful consequence.  

\proclaim{10.2 Lemma} Let $(X,\Delta)$ be a kawamata log terminal pair, where 
$X$ is a projective variety of dimension $n$ and $-(K_X+\Delta)$ is nef.  Let
$\pi:\rmap X.B.$ be the rational map whose existence is guaranteed by (10.1).

 If the MMP holds in dimension $n$ then every component of $\Delta$ dominates $B$.  
\endproclaim
\demo{Proof} Suppose not.  Suppose that $D$ is a component of $\Delta$ which does not dominate 
$B$.  Suppose that we write $\Delta=\Gamma+aD$, where $a$ is the coefficient of $D$ in
$\Delta$.  We run the $(K_X+\Delta)$-MMP with scaling of $D$.  After a series of flips and
divisorial contractions, all of which are $D$-positive, we arrive eventually at a Mori
fibre space $f\:\map W.Z.$.  (Note that even though the pair $(X,\Delta)$ may not
remain kawamata log terminal, nevertheless $K_X+\Gamma$ does remain kawamata log terminal,
and every step of the MMP is a step of the $(K_X+\Gamma)$-MMP).  As every contraction of
this MMP is $D$-positive, $D$ is never contracted and $D$ dominates $Z$.  The fibres
of $f$ are Fano varieties and so in particular they are covered by rational curves.
Pulling these back to $X$, it follows that there is a rational curve $C$ in $X$ which
passes through the general point of $X$ and which intersects $D$.  Hence $C$ is not
contracted under the rational map to $B$.  But then $B$ is covered by rational curves, a
contradiction.  \qed\enddemo

\demo{Proof of (3.9.5)}  Suppose not, suppose that we have an increasing sequence 
$c_i\in \Cal K_n(I)$.  Let $(X_i,\Delta_i)$ be the corresponding sequence of log pairs.
Passing to a log terminal model, we may assume that $X_i$ is $\Bbb Q$-factorial
projective.  By (10.2), restricting to the general fibre of a maximal rationally
chain connected fibration $\pi_i\:\rmap X_i.B_i.$, we may assume that $X_i$ is in fact
rationally connected.

  Suppose that the log discrepancy of $(X_i,\Delta_i)$ is approaching zero.  Then we could
extract a divisor $E_i$ whose coefficient was approaching one.  But this would give us an
increasing sequence in $\Cal N^+_n(I)$, which contradicts (3.9.3).

 Thus we may assume that the log discrepancy of $(X_i,\Delta_i)$ is bounded away from
zero.  As we are assuming (3.8), $\{(X_i,\Delta_i)\}$ forms a bounded family
and the result is clear.  \qed\enddemo

 We also give the following evidence for (3.8):

\proclaim{10.3 Theorem} Assume the MMP in dimension $n$ and (3.7), where
we drop the condition that the Picard number is one (that is we assume the full version of
the conjecture of Alexeev-Borisov).

 Then (3.8) holds, where we replace (3.8.3) and (3.8.4) by 
\roster
\item"(3)$'$" $K_X+\Delta$ is numerically trivial. 
\item"(4)$'$" $-K_X$ is big. 
\endroster
\endproclaim

 Note that we only need the numerically trivial case to prove (3.5) and so
replacing (3) by (3)$'$ is reasonable.  

\proclaim{10.4 Lemma} Assume the MMP in dimension $n$ and let $(X,\Delta)$ 
be a kawamata log terminal pair, where $X$ is projective $\Bbb Q$-factorial, 
such that $K_X+\Delta$ is numerically trivial.  

 Then there is a birational map $\rmap X.Y.$, whose inverse does not contract any
divisors, such that $-K_Y$ is ample.
\endproclaim
\demo{Proof} Pick $\epsilon>0$ sufficiently small so that $K_X+(1+\epsilon)\Delta$ is 
kawamata log terminal.  Suppose that $-K_X$ is not nef.  As $-K_X$ is numerically
equivalent to $\Delta$, it follows that $K_X+(1+\epsilon)\Delta$ is not nef.  Thus running
the $(K_X+(1+\epsilon)\Delta)$-MMP we may assume that $-K_Y$ is nef.  The result now 
follows by applying the base point free Theorem to the nef and big divisor $-K_Y$.  \qed\enddemo

 We now show that running this process in reverse is bounded.  

\proclaim{10.5 Lemma} Assume the MMP in dimension $n$.  Suppose we are 
give a bounded family of log pairs $(Y,\Gamma)$, where $Y$ is projective of dimension $n$,
and $K_Y+\Gamma$ is kawamata log terminal and numerically trivial.

 Then the family of all kawamata log terminal pairs $(X,\Delta)$, for which there is a
birational map $\pi\:\rmap X.Y.$, whose inverse does not contract any divisors, where
$\Delta$ is the log pullback of $\Gamma$, is bounded.
\endproclaim
\demo{Proof} We may suppose, using Noetherian induction and passing to the geometric 
generic fibre, that we have a single log pair $(Y,\Gamma)$, and it suffices to prove in 
this case that there are only finitely many such pairs $(X,\Delta)$. 

 As $K_Y+\Gamma$ is kawamata log terminal, there are only finitely many exceptional
valuations of log discrepancy at most one with respect to $K_Y+\Gamma$.  As the only
divisors extracted by $\pi$ are of log discrepancy at most one, it only remains to prove
that the set of all models $\pi\:\rmap X.Y.$, which extract a set of divisors
corresponding to a fixed set of valuations, is finite.

 Suppose that we have two such extractions $\pi\:\rmap X.Y.$ and $\pi'\:\rmap X'.Y.$.  By
assumption $X$ and $X'$ are isomorphic in codimension one, over $Y$.  Pick an ample
divisor $H'$ on $X'$, and let $H$ be the strict transform of $H'$.  Note that $\pi'$ is
the only such model for which $H'$ is ample.  Consider performing a sequence of
$(K_X+\Delta)$-flops, with respect to $H$.  Such a sequence of flops must terminate, as we
are assuming the existence of the MMP.  But termination only occurs when $H$ becomes
ample, and by uniqueness we must have arrived at $\pi'$.  It follows that any two models
are connected by a sequence of $(K_X+\Delta)$-flops, with respect to some effective
divisor.

 As we are assuming the existence of the MMP, there is at least one model, for which $\pi$
is a morphism.  Pick an ample divisor $H$ on $Y$.  Then we may decompose $\pi^*H$ as
$A+E$, where $A$ is ample and $E$ is effective.  Now $K_X+\Delta+\epsilon E$ is kawamata
log terminal, for $\epsilon$ sufficiently small and it follows that $K_X+\Delta+\epsilon
E$ is relatively ample, for $\epsilon$ positive.  As we assuming the MMP in dimension $n$,
it follows, by the results of either \cite{HK00} or \cite{Shokurov96}, that there is a
decomposition of the relative effective cone of divisors into finitely many chambers, such
that the end result of running a sequence of flops only depends on the chamber to which
the effective divisor belongs.  In particular there are only finitely many possible
models, as required.  \qed\enddemo

\demo{Proof of (10.3)} Applying (10.4) we may assume that there is a 
birational map $\rmap X.Y.$ such that $-K_Y$ is ample.  The log discrepancy of $K_Y$ is at
least the log discrepancy of $K_Y+\Gamma$, where $\Gamma$ is the pushforward of $\Delta$,
which is equal to the log discrepancy of $K_X+\Delta$, which is bounded away from zero, by
assumption.  As we are assuming the strong version of (3.7), which places no
restriction on the Picard number, it follows that we may assume that $(Y,\Gamma)$ forms a
bounded family.  Now apply (10.5). \qed\enddemo

\head \S 11 The toric and toroidal case\endhead 

In this section we prove (3.12).  Most of the proof is completely
straightforward.  All we need to do is go through the proof of (3.9) and check
that we remain in the appropriate category.

 The most technical aspect of the whole proof is the following:

\proclaim{11.1 Lemma} Let $\pi\:\map X.Y.$ be a toroidal and projective morphism, 
where we work locally about a point $p$ on $Y$.  Suppose that $K_X+\Delta$ is kawamata log
terminal.

 Then we may run the $(K_X+\Delta)$-MMP over $Y$.  
\endproclaim
\demo{Proof} By assumption we may find a toric map $\tilde\pi\:\map\tilde X.\tilde Y.$, 
such that $\pi$ and $\tilde\pi$ are formally analytic isomorphic over the point $p$.  In
particular the relative cone of curves of $X/Y$ is naturally isomorphic to the relative
cone of curves of $\tilde X/\tilde Y$.

 Suppose we are given a ray $R$ on $X$ on which $K_X+\Delta$ is negative.  Let $f\:\map
X.Z.$ be the corresponding contraction.  If $f$ is a fibration, there is nothing to prove,
and if the exceptional locus is a divisor, then we replace $X$ by $Z$ and continue in the
usual way.  The only problem is when $f$ is a small contraction.  In this case we need to
establish the existence of the flip of $f$ and to establish the termination of a sequence
of any such flips.

 Let $\tilde f\:\map \tilde X.\tilde Z.$ be the corresponding contraction for $\tilde X$.
Then $\tilde f$ is toric and the existence and termination of toric flips is known.  As an
immediate consequence we get termination of flips for $X$ over $Y$.  To establish
existence is a little more delicate.  Pick an invariant divisor $\tilde S$ for $\tilde X$
such that $\tilde S$ is also negative on $\tilde R$.  Then we may find an invariant
boundary $\tilde B$ on $\tilde X$ such that $K_{\tilde X}+\tilde S+\tilde B$ is negative
on $\tilde R$.  Let $S$ and $B$ be the corresponding divisors on $X$.  It is proved in
\cite{Shokurov02}, see also \cite{Corti00}, that the existence of the flip of $f$ is
equivalent to finite generation of the sheaf $f_*\Cal F$, where $\Cal F$ is a sheaf with
support on $S$.  As the flip of $\tilde f$ exists, it follows that the corresponding sheaf
is finitely generated on $\tilde Z$.  Since $f_*\Cal F$ is supported on the exceptional
locus of $\map Z.Y.$, it follows that $f_*\Cal F$ is also finitely generated, whence the
existence of the flip of $f$ follows from the existence of the flip of $\tilde f$.
\qed\enddemo

 Note that (11.1) is easy to prove over $\Bbb C$, if we define toroidal in
this case to mean that there is an analytic isomorphism, rather than just a formal
analytic isomorphism.

We will also need the following modified version of (3.8).  

\proclaim{11.2 Proposition} Fix an integer $n$ and a positive real 
number $\epsilon$.  

 Then the family of all pairs $(X,\rup \Delta.)$ such that 
\roster 
\item the dimension of $X$ is $n$, 
\item the log discrepancy of the pair $(X,\Delta)$ is at least $\epsilon$, and the 
coefficients of $\Delta$ are at least $\epsilon$,
\item $-(K_X+\Delta)$ is nef, 
\item $X$ is a toric variety, 
\endroster 
is bounded. 
\endproclaim
\demo{Proof} Let $(X,\Delta)$ be any such log pair.  Passing to a log terminal model, we 
may assume that $X$ is a toric $\Bbb Q$-factorial projective variety.  As the invariant
divisors span the Picard group there is an ample divisor which is an effective linear
combination of the invariant divisors.  In particular as $-K_X$ is equivalent to the sum of
the invariant divisors, $-K_X$ is automatically big.  On the other hand any nef divisor on
a toric variety is semi-ample.  Thus, possibly enlarging $\Delta$, we may assume that
$K_X+\Delta$ is numerically trivial.  

 Since the toric MMP is known to exist, we may apply the argument appearing in the proof
of (10.3).  \qed\enddemo

\demo{Proof of (3.12)} ${}^t$(3.5.1) follows just as in the proof of  
(3.9.1), as we may run the MMP in the toric category.  

 To prove ${}^t$(3.5.2), we follow the proof of (3.9.2). By assumption
we may find a pair $(Y,\Gamma)$ which satisfies condition ${}^tT$, such that $c$ is
associated to the pair $(Y,\Gamma)$.  In particular, $(Y,\Gamma)$ has a toroidal
resolution $\pi\:\map W.Y.$.  As condition $T$ is local, we may assume that $Y$ is affine.
Cutting by hyperplanes, we may assume that the only log canonical centre is a point.  Pick
a divisor $X$ in $W$ of log discrepancy zero.  By assumption, $X$ is isomorphic to a toric
variety.  By (11.1) we may run an appropriate MMP and reduce to the case that
$\pi$ extracts only the divisor $X$.  Restricting to $X$ and applying
(4.3), we get a log pair $(X,\Delta)$ where $c$ is associated to $\Delta$
and $X$ is toric. The condition that the pair $(Y,\Gamma)$ has a toroidal resolution
implies that every exceptional valuation of log discrepancy zero is toric.

 To prove ${}^t$(3.5.3-4), note that (3.7) is known for toric
varieties, see \cite{BB92}.  As in the proof of (3.9.3-4), we reduce to the case
where there is a component whose coefficient is approaching one.  In the toric case, we
may further assume that this component is an invariant divisor, so that when we apply
adjunction we do not leave the toric category.

 ${}^t$(3.5.5) follows as in the proof of (3.9.5), using
(11.2). \qed\enddemo

\head Bibliography \endhead
\bibliography{cyclic}

\enddocument